\pdfoutput=1

\documentclass{article}
\usepackage[english]{babel}

\usepackage[letterpaper,top=2cm,bottom=2cm,left=3cm,right=3cm,marginparwidth=1.75cm]{geometry}
\usepackage{algorithm}
\usepackage{algpseudocode}
\usepackage{amsmath}
\usepackage{amssymb}
\usepackage{amsthm}
\usepackage{bbm}
\usepackage{graphicx}
\usepackage[colorlinks=true, allcolors=blue]{hyperref}
\usepackage{indentfirst}
\usepackage{todonotes}
\usepackage[style=numeric-comp,sorting=none]{biblatex}
\theoremstyle{definition}

\newtheorem{property}{Property}[section]

\title{Mixed-Integer Programming decomposition for stochastic programming with quantiles}

\author{Gabriel Gouvine} 

\begin{document}
\maketitle

\begin{abstract}

In this paper, we present a novel approach for solving a maintenance scheduling problem from the French electricity grid company RTE, as presented in the ROADEF 2020 challenge. 

Our approach combines constraint generation with a new family of cutting planes to address the nonconvexity of the problem.
We provide mathematical proofs and separation algorithms for the cutting planes and study the practical impact of our approach on the challenge instances.

This method is applicable to a wide range of problems involving the modeling of uncertainty, and is shown to improve results significantly on real-world instances.

\end{abstract}

\section{Introduction}
The ROADEF 2020 challenge presents a problem arising from RTE's maintenance planning~(\cite{manuel_ruiz_roadef_2020}).
RTE is the company that operates the french electric power grid: reliability is crucial, and any maintenance operation must take into account the risk of the operation for the grid's integrity.
The challenge problem is to find a maintenance schedule taking into account the predicted risk for the power grid.


The uncertainty affecting the grid operations and the energy industry is modeled as a stochastic optimization problem. The optimization objective is a function of the risk's probability distribution, which in practice is approximated over a large but finite number of scenarios.
For this challenge, RTE considers two quality metrics across these scenarios: the mean risk and, in order to penalize worst-case scenarios, a quantile of the risk.
The use of a quantile in the objective function transforms what would be a simple scheduling problem into a complex problem with a non-convex objective.

A useful approach to solve such a two-stage stochastic problem with a mixed-integer programming (MILP) solver is the application of a constraint generation method i.e. Bender's decomposition~(\cite{ruszczynski_decomposition_1997}).
The classical Bender's decomposition is limited to purely continuous linear problems, and a large body of the literature has been dedicated to its generalizations to integer problems~(\cite{fakhri_benders_2017}).

In this paper, after restating the challenge problem in \ref{sec:problem-statement}, we present the MILP model for this problem and the quantile function in \ref{sec:milp-formulation}.
In order to tighten the model, we then introduce a family of valid inequalities for the quantile function in \ref{sec:polyhedral}, as well as an efficient separation algorithm.
To further improve solution time, we design a constraint generation approach in \ref{sec:constraint-generation}, which yields a second family of valid inequalities.
Finally, we test these various methods in practice on the challenge instances in \ref{sec:numerical-results} and provide numerical results, showing that they improve solution time and bounds significantly over the standard MILP model.

The use of a quantile function is common in stochastic programming. In the field of constraint programming, it is known as the max\_n constraint~(\cite{beldiceanu_pruning_2001}), and can be used to implement the more ubiquitous sort constraint~(\cite{older_getting_1995}). We hope that our improvements in its modeling can impact the resolution of many other problems.

\section{Problem statement}
\label{sec:problem-statement}

The goal is to schedule interventions $i \in \mathcal{I}$ according to their starting date $t \in \mathcal{T}$, with a granularity of one day.
Once started, an intervention cannot be preempted.

Personnel and machines are represented as resources $c \in \mathcal C$, and the amount used by the interventions at each timestep must respect the available resource budget.
We denote $\eta_{it'}^{ct}$ the amount of resource $c$ used by intervention $i$ at date $t$ if it starts at $t'$, and $l_t^c$ and $u_t^c$ the minimum and maximum allowed usage for resource $c$ at date $t$.

The planners may define additional requirements, such as forbidding some starting dates for specific interventions, or forbidding some interventions to overlap.
From a theoretical standpoint, these additional constraints are equivalent to resource constraints, and we will conflate them in our notation.

To model the stochastic nature of the problem, the planners define a finite number of scenarios $s \in \mathcal{S}$, and a risk is computed for each timestep and scenario, as the sum of the risks contributed by each intervention.
We denote $\rho_{it'}^{st}$ the risk contributed by intervention $i$ for scenario $s$ at date $t$ if it starts at $t'$.

The objective function is computed separately for each timestep.
To promote low expected risk, it takes into account the mean of the risk across scenarios.
In order to introduce more flexibility, the planners introduce a parameter $k$, and the objective function accounts for the $k^\text{th}$ largest risk.
This allows the planners to penalize bad scenarios while being robust with respect to the worst case. As special cases, they can model the minimum ($k = \lvert \mathcal{S} \rvert$), the median ($k =  \frac{1}{2}\lvert \mathcal{S} \rvert$) and the largest risk ($k = 1$). Since the goal is to penalize worst case scenarios, usually $k \ll \lvert \mathcal{S} \rvert$.

In practice, most instances available from RTE aim at optimizing the 5\% worst case, with $k = \left \lceil 0.05 \lvert \mathcal{S} \rvert \right\rceil$. The number of simulated scenarios $\lvert \mathcal{S} \rvert$ ranges up to 600, with up to 1000 interventions ($\lvert \mathcal{I} \rvert \leq 1000$) to schedule over one year ($\lvert \mathcal{T} \rvert = 365$).

\section{Integer programming model}
\label{sec:milp-formulation}

We model the problem as a mixed-integer programming problem as follows. We denote $Q_k(x), x \in \mathbb{R}^n$ the quantile function, returning the $k^{\text{th}}$ largest component of its argument.

\begin{align}
&\text{minimize } &&\alpha \sum_{t \in \mathcal{T}} m_t + \beta \sum_{t \in \mathcal{T}} q_t \label{model-objective}\\
&\text{subject to } &&\sum_{t \in \mathcal{T}} x_{it} = 1, &i \in \mathcal{I} \label{model-unique-timestep}\\
&&&l_t^c \leq \sum_{i \in \mathcal{I}} \sum_{t' \in \mathcal{T}} \eta_{it'}^{ct} x_{it'} \leq u_t^c, &c \in \mathcal{C}, t \in \mathcal{T} \label{model-resource-constraints}\\
&&&r_{st} = \sum_{i \in \mathcal{I}} \sum_{t' \in \mathcal{T}} \rho_{it'}^{st} x_{it'}, &s \in \mathcal{S}, t \in \mathcal{T} \label{model-define-risk}\\
&&&m_t = \frac{1}{\lvert \mathcal{S} \rvert} \sum_{s \in \mathcal{S}} r_{st}, &t \in \mathcal{T} \label{model-define-mean-risk}\\
&&&q_t \geq Q_k(r_{st}, s \in \mathcal{S}), & t \in \mathcal{T} \label{model-define-quantile-risk}\\
&&&x_{it} \in \{0, 1\}, & i \in \mathcal{I}, t \in \mathcal{T} \label{model-decisions-1}\\
&&&m_t, q_t, r_{st} \in \mathbb{R}, & s \in \mathcal{S}, t \in \mathcal{T} \label{model-decisions-2}
\end{align}

The decision variables $x_{it}$ (\ref{model-decisions-1}) model the start of an intervention at a given timestep.
Decision variables $r_{st}$, $m_t$, $q_t$ (\ref{model-decisions-2}) represent the risk per scenario, the mean risk and the quantile of the risk, respectively.

Constraint \ref{model-unique-timestep} ensures that each intervention is scheduled exactly once.
Constraint \ref{model-resource-constraints} defines the resource constraints.
Constraint \ref{model-define-risk} defines the risk per scenario. Constraint \ref{model-define-mean-risk} defines the mean risk, while constraint \ref{model-define-quantile-risk} 
defines the quantile of the risk; being a minimization objective, an inequality is sufficient.

Constraint \ref{model-define-quantile-risk} is non-convex, as we will prove later, but can be modeled easily.
It states that $q_t$ is larger than the $k^{\text{th}}$ largest element of $\{r_{st}, s \in \mathcal{S}\}$
An equivalent formulation is that there exists a subset $\mathcal{V}_t \subset \mathcal{S}$ of size $\lvert \mathcal{V}_t \rvert = \lvert \mathcal{S} \rvert - k + 1$ such that for all $s \in \mathcal{V}_t, q_t \geq r_{st}$.
We can then add binary decision variables $v_{st}$ to the model such that $v_{st} \Rightarrow s \in \mathcal{V}_t$.
Using an indicator constraint (\ref{model-indicator-constraint}) as supported by many solvers, we use the following well-known model for \ref{model-define-quantile-risk} (\cite{norkin_mixed_2010}):

\begin{align}
&v_{st} \implies q_t \geq r_{st}, & s \in \mathcal{S}, t \in \mathcal{T} \label{model-indicator-constraint}\\
&\sum_{s \in \mathcal{S}} v_{st} = \lvert \mathcal{S} \rvert - k + 1, & t \in \mathcal{T} \\
&v_{st} \in \{0, 1\}, & s \in \mathcal{S}, t \in \mathcal{T}
\end{align}

For solvers that do not support indicator constraints, we need bounds $l_{st} \leq r_{st} \leq u_{st}$ on $r_{st}$ to express a fully linear model using the big-M method.
Here we can pick the naive bounds $l_{st} = 0$ and $u_{st} = \sum_{i \in \mathcal{I}} \sum_{t' \in \mathcal{T}} \rho_{it'}^{st}$.
With the big-M constraint \ref{model-bigm-constraint}, this yields the following integer linear programming model:

\begin{align}
&q_t \geq r_{st} - (u_{st} - l_{st}) (1 - v_{st}), & s \in \mathcal{S}, t \in \mathcal{T} \label{model-bigm-constraint}\\
&\sum_{s \in \mathcal{S}} v_{st} = \lvert \mathcal{S} \rvert - k + 1, & t \in \mathcal{T}\\
&v_{st} \in \{0, 1\}, & s \in \mathcal{S}, t \in \mathcal{T}
\end{align}

\section{Polyhedral analysis}
\label{sec:polyhedral}

The model defined above for the quantile function has a bad linear relaxation.
This makes it difficult for the solver to prove good lower bounds on the objective function or to make accurate decisions during branch-and-bound.
In order to improve the relaxation, we look for valid convex inequalities to be added to the model.

We denote $Q_k(x), x \in \mathbb{R}^n$ the function that returns the $k^\text{th}$ largest component of $x$.
We look for valid convex inequalities for the problem $y \geq Q_k(x)$.

$Q_k$ itself is not convex if $k \ne 1$, since we can find points whose value is lower than the value of their average. For example
$$
Q_k(\overbrace{\ldots, 2}^{k-2}, 1, 1, 0, \ldots) = 1
> Q_k(\overbrace{\ldots, 2}^{k-2}, 2, 0, \ldots) = Q_k(\overbrace{\ldots, 2}^{k-2}, 0, 2, \ldots) = 0
$$

For this reason, no valid convex inequality can be tight. Even more, no convex inequality can be valid for all $x \in \mathbb{R}^n$, or for all $x \in \mathbb{R}_+^n$ except the trivial inequality $y \geq 0$.
To see this, we can write any point in $\mathbb{R}_+^n$ as a convex combination of $n$ points on the coordinate axes, for which all components except one are $0$. $Q_k$ is $0$ on such points if $k \ne 1$, so that $\text{cv }Q_k \leq 0$. Similarly, we can write any point of $\mathbb{R}^n$ as a convex combination of points with arbitrarily low images.
This forces us to look for inequalities defined on a bounded domain instead of general inequalities on $\mathbb{R}^n$.

\begin{property}
Given $x \in \mathbb{R}^n$, $0 \leq x_i \leq u_i$, for any subset $\mathcal{P}$ of $[1\mathrel{{.}\,{.}}n]$ with at least $k$ elements and $U \in \mathbb{R}$, the following inequality is valid for $y \geq Q_k(x)$:
\begin{equation}
\label{eqn-convex-box}
y \geq \frac{U}{\lvert \mathcal{P} \rvert - k + 1} \left(\sum_{i \in \mathcal{P}} \frac{x_i}{\max(U, u_i)} - k + 1\right)
\end{equation}
\end{property}

Restricting the domain of $x$ to the unit cube, we find the more intuitive family of linear inequalities that follows:

\begin{equation}
\label{eqn-convex-unit}
y \geq \frac{1}{\lvert \mathcal{P} \rvert - k + 1}\left(\sum_{i \in \mathcal{P}}x_i - k + 1\right)
\end{equation}

\begin{proof}
As the quantile function is non-decreasing, we only need study the case where $x_i = 0$ for $i \notin \mathcal{P}$.
We assume without loss of generality that $\mathcal{P} = [1 \mathrel{{.}\,{.}} p]$ and $x_i \geq x_{i+1}$ for $1 \leq i < p$.
\begin{align*}  
Q_k(x) = & x_k \\
= & \frac{U}{p - k + 1}\left(\sum_{i=k}^p \frac{x_k}{U}\right) \\
\geq & \frac{U}{p - k + 1}\left(\sum_{i=1}^{k-1} 1 + \sum_{i=k}^p \frac{x_k}{U} - k + 1\right) \\
\geq & \frac{U}{p - k + 1}\left(\sum_{i=1}^p \frac{x_i}{\max(U, u_i)} - k + 1\right) \\
\end{align*}
\end{proof}

This family of inequalities is dominated by the cases where $U = u_i$, which makes it finite but still of exponential size in general.
For this reason, we look for a fast separation algorithm to find the most violated constraint for a given value of $x$.
It can then be used to generate cutting planes during the solution process.

For a given value of $U$, the most violated constraint can be found in linearithmic time ($\mathcal{O}(n \log n)$) by sorting the coordinates of $x$ and iterating on the size of $\lvert \mathcal{P} \rvert$, as implemented by Algorithm 
\ref{alg:symmetric_separation}.

\begin{algorithm}
\caption{Constraint separation from a sorted vector}
\label{alg:symmetric_separation}
\begin{algorithmic}

\Require $\frac{x_i}{\max(U, u_i)} \geq \frac{x_{i+1}}{\max(U, u_{i+1})} \forall i, 1 \leq i < n$

\Function {Constraint Separation}{$U$, $x \in \mathbb{R}^n$}

\State $c \gets \sum_{i=1}^{k-1} \frac{x_i}{\max(U, u_i)}$
\State best bound $\gets -\infty$
\State best p $\gets k$
\For{$p \in [k \mathrel{{.}\,{.}} n]$}
\State $c \gets c + \frac{x_p}{\max(U, u_p)}$ \Comment Maintain $\sum_{i=1}^p\frac{x_i}{\max(U, u_i)}$
\State bound $\gets \frac{c - k + 1}{p - k + 1}$  \Comment Best lower bound for $\lvert \mathcal{P} \rvert = p$
\If{bound $>$ best bound}
\State best bound $\gets$ bound
\State best p $\gets p$
\EndIf
\EndFor

\State \Return $[1 \mathrel{{.}\,{.}} \text{best p}]$ \Comment $\mathcal{P}$ with the best bound

\EndFunction
\end{algorithmic}
\end{algorithm}

In the general case, we can repeat this algorithm for all possible values of $U$, sorting by $\frac{x_i}{\max(U, u_i)}$.
Once the vector is sorted, the separation algorithm runs in linear time for each $U$.
If we perform two initial sorts, by $x_i$ and $\frac{x_i}{u_i}$, they can be used to sort the vector by $\frac{x_i}{\max(U, u_i)}$ in linear time using one merge-sort step.
This gives a constraint separation algorithm that runs in quadratic time ($\mathcal{O}(n^2)$) in the general case.

\section{Constraint generation method}
\label{sec:constraint-generation}

Another limitation of the mixed-integer linear formulation is its size. The model of the quantile function introduces one binary variable and one indicator constraint for each element of the population.
For the problem at hand, this formulation is repeated for each timestep, making the model grow even faster.
This repeated subproblem is a natural target for decomposition methods.

We propose a constraint generation method with subproblems of the form $y \geq Q_k(Ax)$, where $A \in \mathbb{R}_+^{n \times m}$ and $x \in \{0, 1\}^m$. In the problem at hand, $x$ maps to the decision variables $x_{it}$, and the entries of $A$ are the parameters $\rho_{it'}^{st}$ defining the risk per scenario.

In its simplest form, we add a new constraint for each incumbent solution $\widetilde{x}$, which enforces the correct value of $y$ at this point.
As the set of possible values for $x$ is finite, this guarantees eventual convergence.
This is similar to the addition of a ``no-good'' cut in \cite{kim_solving_2002} and \cite{codato_combinatorial_2006}, with the addition of the continuous variable $y$.
One simple such constraint is:

\begin{equation}
\label{eqn-nogood-constraint}
y \geq Q_k(A\widetilde{x}) \left(\sum_{j=1}^m \widetilde{x_{j}} x_j - \sum_{j=1}^m \widetilde{x_j} + 1\right)
\end{equation}

This constraint enforces $y \geq Q_k(A\widetilde{x})$ for $x \geq \widetilde{x}$, which is correct since $x \mapsto Q_k(Ax)$ is non-decreasing.
With large coefficients and a huge negative constant term, this first constraint is very weak and leads to no improvements in the linear relaxation in most cases.
Moreover, if there is a $j$ such that $x_{j} = 0 $ and $\widetilde{x_{j}} = 1$, the right-hand side is always non-positive: the constraint is only active on solutions $x \geq \widetilde{x}$.
Finally, its right-hand side is never larger than $Q_k(A\widetilde{x})$.
In order to obtain tighter bounds and a faster solution process, we introduce the following family of valid inequalities.

\begin{property}
Given $A \in \mathbb{R}^{n \times m}$ and $x \in \mathbb{R}^n$, $l_i \leq x_i \leq u_i$, for any subset $\mathcal{P}$ of $[1\mathrel{{.}\,{.}}n]$ with at least $k$ elements and $\beta \in [0,1]^n$, the following inequality is valid for $y \geq Q_k(Ax)$:

\begin{align}
\begin{aligned}
\label{eqn-subset-constraint}
y \geq& 
\min_{i \in \mathcal{P}} \sum_{j=1}^m a_{ij} \left((1 - \beta_j) l_j +\beta_j u_j\right) \\
&+ \sum_{j=1}^m \min_{i \in \mathcal{P}} (a_{ij}) (1 - \beta_j) (x_j - l_j) \\
&+ \sum_{j=1}^m \max_{i \in \mathcal{P}}(a_{ij}) \beta_j (x_{j} - u_j) \\
\end{aligned}
\end{align}
\end{property}

\begin{proof}
\begin{align*}
Q_k(Ax) \geq& \min_{i \in \mathcal{P}} \sum_{j=1}^m a_{ij} x_{j} \\
Q_k(Ax) \geq& \min_{i \in \mathcal{P}} \sum_{j=1}^m a_{ij} (1 - \beta_j) (x_j - l_j) + a_{ij} \beta_j (x_j - u_j) + a_{ij} \left((1 - \beta_j) l_j +\beta_j u_j\right)\\
Q_k(Ax) \geq& \sum_{j=1}^m \min_{i \in \mathcal{P}} \left( a_{ij} (1 - \beta_j) (x_j - l_j)\right) + \sum_{j=1}^m \min_{i \in \mathcal{P}} \left( a_{ij} \beta_j (x_j - u_j) \right) \\
&+ \min_{i \in \mathcal{P}} \sum_{j=1}^m a_{ij} \left((1 - \beta_j) l_j +\beta_j u_j\right) \\
\intertext{and since $(1-\beta_j)(x_j - l_j) \geq 0$ and $\beta_j(x_j - u_j) \leq 0$:}\\
Q_k(Ax) \geq& \sum_{j=1}^m \min_{i \in \mathcal{P}} (a_{ij}) (1 - \beta_j) (x_j - l_j) + \sum_{j=1}^m \max_{i \in \mathcal{P}}(a_{ij}) \beta_j (x_{j} - u_j) \\
&+ \min_{i \in \mathcal{P}} \sum_{j=1}^m a_{ij} \left((1 - \beta_j) l_j +\beta_j u_j\right)
\end{align*}
\end{proof}

For $A \in \mathbb{R}_+^{n \times m}$ and $x \in \{0, 1\}^m$, we obtain a stronger version of the no-good constraint \ref{eqn-nogood-constraint}. Given an incumbent solution $\widetilde{x}$ and $\mathcal{P}$ such that $Q_k(A\widetilde{x}) = \min_{i \in P} \sum_{j} a_{ij} \widetilde{x_j}$:

\begin{equation}
\label{eqn-strongnogood-constraint}
y \geq Q_k(A\widetilde{x}) +
\sum_{j: \widetilde{x_j} = 0} x_{j} \min_{i \in \mathcal{P}} a_{ij} +
\sum_{j: \widetilde{x_j} = 1} (x_{j}-1) \max_{i \in \mathcal{P}} a_{ij}
\end{equation}

Which, for $\widetilde{x} = 0$, becomes the subset constraint:

\begin{equation}
\label{eqn-simple-subset-constraint}
y \geq \sum_{j=1}^m \min_{i \in \mathcal{P}}(a_{ij}) x_{j}
\end{equation}

Like \ref{eqn-nogood-constraint}, \ref{eqn-strongnogood-constraint} enforces the correct value of $y$ for the incumbent solution $\widetilde{x}$.
However, it provides more information to the solver, with a right-hand side that can become larger than $Q_k(A\widetilde{x})$.
With smaller coefficients for coordinates $j$ where $\widetilde{x_j} = 1$, it can remain useful even when $x \geq \widetilde{x}$ does not hold.

This new family is much more general, and our goal is to use it as cutting planes at the root node of the solution process.
However, given a solution to the linear relaxation, finding the most violated constraint even in the simpler family \ref{eqn-simple-subset-constraint} is NP-hard.

To prove it, consider an undirected graph with $n$ vertices and $m$ edges. We consider the matrix $A$ with $a_{ij} = 0$ if the edge $j$ is incident to node $i$, $1$ otherwise.
$\min_{i \in P}(a_{ij})$ is $1$ if and only if the subset of nodes associated with $P$ does not cover edge $j$.
Finding the largest value for the right-hand side amounts to finding a subset of $k$ vertices that covers the smallest number of edges.
This is the MINIMUM PARTIAL VERTEX COVER problem, which is NP-hard~(\cite{guo_parameterized_2007}).

\section{Numerical results}
\label{sec:numerical-results}

We compare the following resolution methods:
\begin{itemize}
\item The original model with indicator constraints from Section \ref{sec:milp-formulation} (Full)
\item The original model with cutting planes \ref{eqn-convex-box} used to tighten the subproblems $q_t \geq Q_k(r_{st}, s \in \mathcal{S})$ (Full+Sep).
The most violated inequality is added at the root node until no violated inequality remains.
\item Our constraint generation method using \ref{eqn-strongnogood-constraint} (CGen).
\item The two formulations, with the subset constraints \ref{eqn-simple-subset-constraint} (Full+Subs/CGen+Subs) added heuristically at model creation time.
For each intervention and timestep, we use the subset $\mathcal{P}$ corresponding to the $k$ largest risks.
\item In order to compare this heuristic against a best-case scenario, the constraint generation model with the most violated cutting planes \ref{eqn-subset-constraint}
added at the root node (CGen+Best).
We perform 20 rounds of constraint separation at the root node using a MILP subproblem. The time required to solve this subproblem is large and is not included in the measurement.
\end{itemize}

We implement our models in Python and solve them using CPLEX 12.9.
We run our models on the instance sets A, B and C from RTE, for one hour not including model creation time. The code and instances are publicly available (see \ref{sec:statements}).

In the results presented here, we remove instances where all the models reach the optimal solution (8 from A set), as well as instances where all models crash or finish without a feasible solution (4 from B and C sets).
Our benchmarking machine has a 4-core CPU with 10GB of memory.
We report the optimality gap found by the solver in Table \ref{Tab:gap}.
Additional results are provided in Appendix \ref{app:numerical-results}, Table \ref{Tab:lower-bound} and \ref{Tab:solution}, respectively comparing the proven lower bound and the solution found to the best known solution.

\begin{table}[ht]
\begin{center}
\begin{tabular}{ |l|c|c|c|c|c|c| } 
\hline
Instance& Full& Full+Sep& Full+Subs& CGen& CGen+Subs& CGen+Best \\
\hline
A\_02& 54.65\%& 62.64\%& 2.41\%& 2.70\%& \textbf{2.06\%}& \textit{1.79\%} \\
A\_05& 9.27\%& 9.90\%& - & 7.70\%& \textbf{5.99\%}& \textit{5.19\%} \\
A\_08& 6.31\%& 6.98\%& 2.77\%& \textbf{0.00\%}& \textbf{0.00\%}& \textit{0.00\%} \\
A\_11& 8.72\%& 10.17\%& 8.63\%& \textbf{4.35\%}& 4.63\%& \textit{4.16\%} \\
A\_13& \textbf{0.10}\%& \textbf{0.10}\%& 0.12\%& \textbf{0.10}\%& 0.12\%& \textit{0.08\%} \\
A\_14& 9.15\%& 9.82\%& 9.07\%& 6.72\%& \textbf{6.44}\%& \textit{6.17\%} \\
A\_15& 12.78\%& 14.42\%& 12.10\%& \textbf{7.37}\%& 7.53\%& \textit{6.44\%} \\
\hline
B\_01& - & - & 8.24\%& 14.09\%& \textbf{7.45}\% &\\
B\_02& \textbf{60.65}\%& - & 62.09\%& 60.79\%& 61.97\% &\\
B\_03& \textbf{68.50}\%& - & - & - & - &\\
B\_04& - & 62.90\%& - & \textbf{9.48}\%& - &\\
B\_05& 51.15\%& 52.86\%& - & 16.62\%& \textbf{5.70\%} &\\
B\_06& 59.86\%& 61.81\%& \textbf{58.98\%}& 60.94\%& 62.05\% &\\
B\_07& 47.45\%& - & - & \textbf{46.31\%}& 46.60\% &\\
B\_08& 68.08\%& - & \textbf{4.02\%}& 7.92\%& 6.12\% &\\
B\_09& 67.44\%& - & \textbf{3.46\%}& 22.37\%& 4.59\% &\\
B\_10& - & - & - & \textbf{47.81\%}& - &\\
B\_11& 58.58\%& - & \textbf{56.05\%}& 56.57\%& 57.61\% &\\
B\_12& - & - & - & \textbf{12.40\%}& - &\\
B\_13& - & - & - & 12.62\%& \textbf{7.13\%} &\\
B\_14& - & - & - & \textbf{39.17\%}& - &\\
\hline
C\_01& 70.65\%& 68.46\%& \textbf{5.47\%}& 77.31\%& 6.45\% &\\
C\_02& - & - & - & \textbf{39.34\%}& 39.63\% &\\
C\_03& 70.03\%& \textbf{67.23\%}& - & 71.80\%& - &\\
C\_04& 64.85\%& 64.29\%& - & \textbf{9.38\%}& - &\\
C\_05& 59.77\%& 60.70\%& 4.30\%& 14.44\%& \textbf{4.20\%} &\\
C\_06& 54.12\%& - & - & \textbf{47.55\%}& 47.61\% &\\
C\_07& 54.23\%& - & \textbf{3.77\%}& 16.50\%& 6.06\% &\\
C\_08& 50.28\%& - & - & \textbf{44.77\%}& 44.80\% &\\
C\_09& \textbf{69.92\%}& - & 70.30\%& 70.03\%& 71.10\% &\\
C\_10& - & - & - & \textbf{70.37\%}& - &\\
C\_11& 59.58\%& - & 9.32\%& 12.64\%& \textbf{7.03\%} &\\
C\_12& - & - & - & \textbf{47.35\%}& - &\\
\hline
\end{tabular}
\end{center}
\caption{Optimality gap obtained by the solver with various methods after one hour.}
\label{Tab:gap}
\end{table}

The quality of the constraints that can be deduced from the convex relaxation of the quantile function alone is highly dependent on the tightness of the lower and upper bounds on the risk values.
In our experiments, the convex relaxation of the quantile function leads to an actually worse model than the natural MILP model, presumably due to the time spent in the constraint separation algorithm and the computational cost of the additional constraints.
This is expected in hindsight, as the interventions tend to be spread over time due to the resource constraints, and the risk values $r_{st}$ are much lower than the theoretical upper bound, where the inequality \ref{eqn-convex-box} is tighter.

Looking at the objective value only, the constraint generation method is better than the full model for 31 of 33 instances. For 19 instances, one of the two constraint generation methods obtains a solution within 1\% of the best known solution for the challenge, against only 7 for the full model. The solution of instance A\_08 is even proven optimal.

The effect on the proven lower bound is less clear: while the lower bound is better for all small A instances, the full model with our subset constraints often proves tighter lower bounds on B and C instances.

The effect of adding subset constraints at the root node varies between instances. For some of them, the gap between the lower bound and the best known solution is reduced by a factor of 20, while for others it is virtually unchanged. It is particularly efficient when combined with the full model, where it yields some of the best bounds.

We see that the ``ideal'' case, with perfect constraint separation at the root node, is still significantly better at reducing the optimality gap. With appropriate fast heuristics to find the cutting planes, there is certainly room for improvement.
Another consequence of our crude heuristic is that the model is much larger. The solver runs out of memory or is unable to find a solution on our machine for some instances, making the pure constraint generation method a credible alternative if finding solutions is more important than the quality of the bounds.

\section{Further work}

Our work on this approach doesn’t end here, as neither the constraint generation algorithm nor the cutting planes are restricted to the challenge problem.
The quantile function is used extensively in stochastic programming, and our method could be applied to many other problems in the field.

In particular, the quantile function we studied can be used to implement and strengthen the sort constraint in constraint programming. Indeed, for $x, y \in \mathbb{R}^n$, $y = \text{sort}(x)$ can be rewritten as:
\begin{align*}
    y_k &\geq Q_k(x), & 1 \leq k \leq n \\
    -y_k &\leq Q_{n-k+1}(-x), & 1 \leq k \leq n 
\end{align*}
Therefore, all valid linear inequalities above can be used on linear programming reformulations of the sort constraint. Other work could apply them to these other classes of problems, or integrate them in reformulation library.

Finally, the more general family of inequalities presented here for $y \geq Q_k(Ax)$ (\ref{eqn-subset-constraint}) can be applied to $y \geq Q_k(x)$ (\ref{eqn-convex-box}) by setting $A$ to the identity matrix.
Some of the inequalities can be derived with both formulas. For example, when $u = 1$ and $l = 0$, they achieve the same result for $\lvert \mathcal{P} \rvert = k$ and $\beta = 1$. However, \ref{eqn-convex-box} appears to be weaker for $\lvert \mathcal{P} \rvert > k$. While we were not able to find a generalization, this hints that even stronger general
inequalities are possible.

\section{Conclusion}

In this paper, we presented a mathematical programming method for a stochastic programming problem. The modeling of uncertainty using a quantile function makes the problem hard to model, and the natural mixed integer programming model is large and slow to converge. We followed two paths to improve its modeling.
First, we studied valid inequalities for the quantile function and developed separation algorithms and heuristics; this led to a stronger linear relaxation of the model. Second, we introduced a constraint generation approach; this method made the model much smaller, making the solution process faster and tackling larger instances.

We then studied the practical performance of our methods on the challenge instances. The experiments show that our constraint generation approach obtains much better solutions, and that the new cutting planes obtain better optimality bounds compared to the natural model.

\printbibliography

\clearpage
\appendix
\section{Numerical results}
\label{app:numerical-results}

\begin{table}[ht]
\begin{center}
\begin{tabular}{ |l|c|c|c|c|c| }
\hline
Instance& Full& Full+Subs& CGen& CGen+Subs& CGen+Best\\
\hline
A\_02& 54.47\%& 2.05\%& 2.38\%& \textbf{1.91\%}& \textit{1.61\%}\\
A\_05& 6.76\%& - & 6.22\%& \textbf{5.43\%}& \textit{4.33\%}\\
A\_08& 6.25\%& 2.74\%& \textbf{0.00\%}& \textbf{0.00\%}& \textit{0.00\%}\\
A\_11& 6.33\%& 5.89\%& \textbf{4.18\%}& 4.45\%& \textit{3.88\%}\\
A\_13& 0.08\%& 0.09\%& \textbf{0.05\%}& 0.06\%& \textit{0.05\%}\\
A\_14& 7.44\%& 7.06\%& \textbf{5.55\%}& 5.76\%& \textit{4.91\%}\\
A\_15& 7.76\%& 7.60\%& \textbf{6.66\%}& 6.74\%& \textit{5.19\%}\\
\hline
B\_01& - & \textbf{6.22\%}& 13.22\%& 7.28\%& \\
B\_02& 50.24\%& \textbf{50.11\%}& 55.12\%& 55.13\%& \\
B\_03& \textbf{63.58\%}& - & - & - & \\
B\_04& - & - & \textbf{9.34\%}& - & \\
B\_05& 47.12\%& - & 15.24\%& \textbf{5.48\%}& \\
B\_06& \textbf{50.04\%}& 50.13\%& 55.04\%& 55.09\%& \\
B\_07& \textbf{41.50\%}& - & 41.58\%& 41.63\%& \\
B\_08& 64.09\%& \textbf{3.88\%}& 7.90\%& 6.11\%& \\
B\_09& 66.03\%& \textbf{3.24\%}& 15.73\%& 4.59\%& \\
B\_10& - & - & \textbf{46.04\%}& - & \\
B\_11& \textbf{48.25\%}& 49.65\%& 52.51\%& 52.45\%& \\
B\_12& - & - & \textbf{11.18\%}& - & \\
B\_13& - & - & 12.25\%& \textbf{7.12\%}& \\
B\_14& - & - & \textbf{37.61\%}& - & \\
\hline
C\_01& 66.61\%& \textbf{4.73\%}& 76.44\%& 6.40\%& \\
C\_02& - & - & \textbf{37.75\%}& 37.86\%& \\
C\_03& \textbf{65.14\%}& - & 66.57\%& - & \\
C\_04& 63.97\%& - & \textbf{9.29\%}& - & \\
C\_05& 57.35\%& \textbf{2.68\%}& 11.70\%& 3.83\%& \\
C\_06& \textbf{44.85\%}& - & 45.31\%& 45.39\%& \\
C\_07& 51.70\%& \textbf{3.45\%}& 15.06\%& 5.95\%& \\
C\_08& \textbf{41.45\%}& - & 42.66\%& 42.65\%& \\
C\_09& 59.58\%& \textbf{59.53\%}& 65.24\%& 65.23\%& \\
C\_10& - & - & \textbf{70.18\%}& - & \\
C\_11& 54.45\%& 7.42\%& 11.96\%& \textbf{7.00\%}& \\
C\_12& - & - & \textbf{45.51\%}& - & \\
\hline
\end{tabular}
\end{center}
\caption{Gap between the best lower bound proven by the solver in one hour and the best solution reported for the challenge}
\label{Tab:lower-bound}
\end{table}

\begin{table}[ht]
\begin{center}
\begin{tabular}{ |l|c|c|c|c|c| }
\hline
Instance& Full& Full+Subs& CGen& CGen+Subs& CGen+Best\\
\hline
A\_02& 0.39\%& 0.38\%& 0.33\%& \textbf{0.16\%}& \textit{0.18\%}\\
A\_05& 2.76\%& -& 1.60\%& \textbf{0.60\%}& \textit{0.91\%}\\
A\_08& 0.06\%& 0.03\%& \textbf{0.00\%}& \textbf{0.00\%}& \textit{0.00\%}\\
A\_11& 2.61\%& 3.00\%& \textbf{0.18\%}& 0.19\%& \textit{0.30\%}\\
A\_13& \textbf{0.02\%}& \textbf{0.02\%}& 0.05\%& 0.06\%& \textit{0.03\%}\\
A\_14& 1.88\%& 2.21\%& 1.25\%& \textbf{0.72\%}& \textit{1.35\%}\\
A\_15& 5.74\%& 5.12\%& \textbf{0.76\%}& 0.86\%& \textit{1.33\%}\\
\hline
B\_01& -& 2.20\%& 1.01\%& \textbf{0.18\%}& \\
B\_02& 26.43\%& 31.60\%& \textbf{14.45\%}& 17.98\%& \\
B\_03& \textbf{15.63\%}& -& -& - & \\
B\_04& -& -& \textbf{0.15\%}& - & \\
B\_05& 8.24\%& -& 1.66\%& \textbf{0.24\%}& \\
B\_06& 24.45\%& 21.58\%& \textbf{15.10\%}& 18.36\%& \\
B\_07& 11.32\%& -& \textbf{8.81\%}& 9.30\%& \\
B\_08& 12.52\%& 0.14\%& \textbf{0.01\%}& \textbf{0.01\%}& \\
B\_09& 4.33\%& 0.22\%& 8.56\%& \textbf{0.00\%}& \\
B\_10& -& -& \textbf{3.39\%}& - & \\
B\_11& 24.94\%& 14.56\%& \textbf{9.35\%}& 12.15\%& \\
B\_12& -& -& \textbf{1.39\%}& - & \\
B\_13& -& -& 0.42\%& \textbf{0.01\%}& \\
B\_14& -& -& \textbf{2.58\%}& - & \\
\hline
C\_01& 13.77\%& 0.79\%& 3.86\%& \textbf{0.05\%}& \\
C\_02& -& -& \textbf{2.61\%}& 2.93\%& \\
C\_03& \textbf{16.33\%}& -& 18.55\%& - & \\
C\_04& 2.50\%& -& \textbf{0.10\%}& -&  \\
C\_05& 6.01\%& 1.69\%& 3.21\%& \textbf{0.38\%}& \\
C\_06& 20.21\%& -& 4.29\%& \textbf{4.23\%}& \\
C\_07& 5.54\%& 0.33\%& 1.72\%& \textbf{0.12\%}& \\
C\_08& 17.75\%& -& \textbf{3.83\%}& 3.89\%& \\
C\_09& 34.38\%& 36.28\%& \textbf{16.01\%}& 20.31\%& \\
C\_10& -& -& \textbf{0.66\%}& - & \\
C\_11& 12.69\%& 2.10\%& 0.78\%& \textbf{0.03\%}& \\
C\_12& -& -& \textbf{3.49\%}& - & \\
\hline
\end{tabular}
\end{center}
\caption{Gap between the solution found by the solver in one hour and the best solution reported for the challenge}
\label{Tab:solution}
\end{table}

\end{document}